\documentclass[12pt]{article}

\usepackage{mathtools}
\usepackage{jfe_preambule}
\usepackage{booktabs}
\usepackage{bbm}
\usepackage[]{algorithm2e}
\usepackage{nccmath}
\newcommand{\norm}[1]{\left\lVert#1\right\rVert}

\SetKwInput{kwInit}{Init}
\SetKw{Continue}{continue}
\SetKw{Axis}{axis}

\begin{document}

\setlist{noitemsep}  

\title{Seven proofs of the Pearson Chi-squared independence test and its graphical interpretation}

\author{Eric Benhamou 
\thanks{A.I. SQUARE CONNECT, 35 Boulevard d'Inkermann 92200 Neuilly sur Seine, France and LAMSADE, Université Paris Dauphine, Place du Maréchal de Lattre de Tassigny,75016 Paris, France. E-mail: eric.benhamou@aisquareconnect.com, eric.benhamou@dauphine.eu} 
\textsuperscript{,}
 \textsuperscript{\ddag}
\hspace{0.2cm}, 
Valentin Melot
\thanks{Ecole Normale Supérieure, 45 Rue d'Ulm, 75005 Paris, France. E-mail: valentin.melot@ens.fr}
\textsuperscript{,}
\thanks{the authors would like to mention a fruitful discussion with Robin Ryder that originates this work}
}

\date{}              

\singlespacing

\maketitle

\vspace{-.2in}
\begin{abstract}
This paper revisits the Pearson Chi-squared independence test. 
After presenting the underlying theory with modern notations and showing new way of deriving the proof, we describe
an innovative and intuitive graphical presentation of this test. This enables not 
only interpreting visually the test but also measuring how close or far we are from
accepting or rejecting the null hypothesis of non independence.\end{abstract}

\medskip



\thispagestyle{empty}

\clearpage

\onehalfspacing
\setcounter{footnote}{0}
\renewcommand{\thefootnote}{\arabic{footnote}}
\setcounter{page}{1}


\section{Introduction}
A natural and common question in statistics is to state if two nominal (categorical) variables are independent or not. 
The traditional approach is to use the Pearson Chi-Square test of independence as developed in \cite{Pearson_1900}. We can determine if there is a (statistically) significant relationship between these two nominal variables.  Traditionally, the data are displayed in a contingency table where each row represents a category for one variable and each column represents a category for the other variable.  

For example, say a researcher wants to examine the relationship between gender (male vs. female) and driver riskiness (dangerous or safe driver).  The null (respectively the alternative) hypothesis for this test is that there is no relationship (respectively a relationship) between the two variables: gender and driver riskiness. 

The chi-square test yields only an approximated p-value as this is an asymptotic test as we will see shortly. Hence, this only works when data-sets are large enough. For small sample sizes, Fisher’s (as explained in \cite{Fisher_1922}) or Barnard's (as presented in \cite{Barnard_1945} or in \cite{Barnard_1947}) exact tests are more appropriate but more complex. The large interest of the  Pearson Chi-Square test of independence is its simplicity and robustness as it only relies on two main assumptions: large sample size and independence of observations. It is a mainstream test, available in the core library of \textbf{\textsf{R}}: function \textit{chisq.test} or in python (function \textit{stats.pearsonr} of the scipy library). A natural question is then to graphically represent this test to illustrate if we are close or not to the null hypothesis. Although there has been lot of packages to represent contingency tables, there has been a lack of thought for representing this test visually. This has motivated us to write this work that leads to the writing of a short function in python to do it. We also revisited the proof done by Pearson in 1900 and show that this proof can be derived more elegantly with more recent mathematical tools, namely Cochran theorem that was only found in 1934 (see \cite{cochran_1934}) and also the Sherman Morison matrix inversion formula (provided in 1949 by \cite{Sherman_Morison_1949}). 

The contribution of our paper are twofold: to provide first a modern proof of the Pearson chi square test and second a nice and graphical interpretation of this test. 
Compared to previous proofs as for instance in \cite{Buonocore_Pirozzi_2014}, we are the first one to provide seven proofs for this seminal results with the use of a wide range of tools, like not only Cochran theorem but Sherman Morison formula, Sylvester's theorem and an elementary proof using De Moivre-Laplace theorem. We are also the first paper to suggest to use confidence interval in mosaic plots.

The paper is organized as follows. We first present with modern notation the underlying theory and the seven proofs. We then examine how to display on a single graphic both the contingency tables and the test. We show that this is surprisingly simple yet powerful. We conclude with various examples to illustrate our graphical representation.

\section{Null hypothesis asymptotic}

Let $X_1, X_2, \ldots$ be independent samples from a multinomial(1, $p$) distribution, where $p$ is a $k$-vector
with nonnegative entries that sum to one. That is,
\begin{equation}
P(X_{ij} = 1) = 1 - P(X_{ij} = 0) = p_j   \quad \text{for all}  \> 1 \leq  j \leq k 
\end{equation}
and each $X_i$ consists of exactly $k-1$ zeros and a single one, where the one is in the component
of the “success” category at trial $i$. 

This equation implies in particular that  $ \text{Var}(X_{ij}) = p_j (1-p_j)$. Furthermore, 
$\text{Cov} (X_{ij} ,X_{il}) = \mathbb{E}[ X_{ij} X_{il} ]- p_j p_l = - p_j p_l$ for $ j \neq l$. Therefore, the random vector $X_i$ has covariance matrix given by
\begin{equation}
\Sigma = 
\left( 
\begin{array}{l l l l}
{p_1 (1-p_1)}	& {-p_1 p_2}	& \ldots 		& {-p_1 p_k}		\\
{-p_1 p_2}		& {p_2 (1-p_2)}	& \ldots 		& {-p_2 p_k}		\\
\vdots 			& \vdots 		& \ddots 	& \vdots 		  	\\
{- p_1 p_k } 		& {- p_2 p_k }  	& \ldots 		& {p_k (1-p_k)}
\end{array}
\right) 
\end{equation}

Let us prove shortly that the asymptotic distribution of the Pearson chi-square statistic given by 
\begin{equation}\label{ChiSquareStat}
\chi^2 = \sum _{j=1}^{k} \frac{ (N_j - n p_j)^2 }{ n p_j}
\end{equation}
where $N_j$ is the random variable $n \bar{X}_j$, the number of successes in the $j$th category for trials
$1, \ldots, n$ converges in distribution to the chi-square distribution with $k - 1$ degrees of freedom. This will imply in particular that to test that two samples are from the same statistics, we can use a test of goodness of fit.

The proof is rather elementary and we provide below seven different methods. These proofs show that they are profound connections between binomial, multinomial, Poisson, normal and chi squared distribution for asymptotic cases. They also illustrate that this problem can be tackled with multiple mathematical tools like De Moivre-Laplace theorem that is an early and simpler version of the Central Limit theorem and a recursive induction, but also characteristic function and Lévy's continuity theorem, geometry and linear algebra reasoning that are at the foundation of the Cochran theorem. 

\section{Seven different proofs for the Pearson independence test}
\subsection{First Proof: Sherman Morison formula for direct computation}

Since $\mathbb{E}[X_i] = p$, the central limit theorem implies 
\begin{equation}
\sqrt n ( \bar{X}_n - p ) \xrightarrow[n \rightarrow \infty]{\text{d}} N_k(0, \Sigma)
\end{equation}

where the notation $\xrightarrow[n \rightarrow \infty]{\text{d}}$  
indicates convergence in distribution and $N_k(0,\Sigma)$ is the multi dimensional normal with $k$ dimension and $\Sigma$ as its covariance matrix. Note that $\Sigma$ is not invertible as the sum of any $j$th column of $\Sigma$ is null since it is equal to $p_j - p_j (p_1+ \ldots + p_k)$. 

A first way to tackle the problem of inferring the distribution of the $\chi^2$ statistic is to remove one dimension to the $X$ vector to have a covariance matrix of full rank. More precisely, let us define for each sample $i$, the truncated vector variable $X^{*}_i = (X_{i1}, \ldots  ,X_{i,k-1})^T$. It is the $k-1$ vector consisting of the first $k- 1$ components of $X_i$.  Its covariance matrix is the sub matrix of $\Sigma$ reduced to its first $k-1$ rows and columns. We call it $\Sigma^{*}$. $\Sigma^{*}$ writes as the sum of two simple matrices
\begin{equation}
\Sigma^* = 
\underbrace{\left( 
\begin{array}{l l l l}
{p_1} 	& {0}		& \ldots 		& {0}		\\
{0}		& {p_2 }	& \ldots 		& {0}		\\
\vdots 	& \vdots 	& \ddots 	& \vdots 	\\
{0}		& {0}		& \ldots 		& {p_{k-1}}
\end{array}
\right) }_{A}
-
\underbrace{\left( 
\begin{array}{l l l l}
{p_1 }		\\
{p_2 }		\\
\vdots 		\\
{p_{k-1} } 		
\end{array}
\right) }_{b}
\underbrace{\left( 
\begin{array}{l l l l}
{p_1 }		\\
{p_2 }		\\
\vdots 		\\
{p_{k-1} } 		
\end{array}
\right) ^T}_{b^T}
\end{equation}

We have trivially that $1- b^T A^{-1} b = \sum_{i=1}^k p_i - \sum_{i=1}^{k-1} p_i =  p_k$ and $A^{-1} b = (1, \ldots, 1 ) ^T$ where the last vector is of $k-1$ dimension. $\Sigma^{*}$ inverse, denoted by $(\Sigma^*) ^{-1}$, is therefore given by the Sherman-Morrison formula 
\begin{equation}
(\Sigma^*) ^{-1} = 
\underbrace{\left( 
\begin{array}{l l l l}
{1/p_1} 	& {0}		& \ldots 		& {0}		\\
{0}		& {1/p_2 }	& \ldots 		& {0}		\\
\vdots 	& \vdots 	& \ddots 	& \vdots 	\\
{0}		& {0}		& \ldots 		& {1/p_{k-1}}
\end{array}
\right)}_{A^{-1} }
+
\underbrace{\frac{1}{p_k}  \left(
 \begin{array}{l l l l}
{1}	 	& {1}		& \ldots 		& {1}		\\
{1}		& {1}		& \ldots 		& {1}		\\
\vdots 	& \vdots 	& \ddots 	& \vdots 	\\
{1}		& {1}		& \ldots 		& {1}
\end{array}\right)}_{( 1 / ( 1- b^T A^{-1} b) ( A^{-1} b b^T A^{-1} )}
\end{equation}

Let us write also $p^*$ the $k-1$ vector of probability $p^{*}= ( p_1, \ldots, p_{k-1})^T$. 
The $\chi^2$ statistic of equation (\ref{ChiSquareStat}) can be reformulated as follows:
\begin{eqnarray}
\chi^2 & =& n \sum _{j=1}^{k} \frac{ (\bar X_j - p_j)^2 }{ p_j} \\
& =& n \sum _{j=1}^{k-1} \frac{ (\bar X_j - p_j)^2 }{ p_j} +  \frac{ (\bar X_k - p_k)^2 }{ p_k} \\
& =& n \sum _{j=1}^{k-1} \frac{ (\bar X_j - p_j)^2 }{ p_j} +  \frac{ (\sum _{j=1}^{k-1} (\bar X_j - p_j) )^2 }{ p_k} 
\end{eqnarray}

where we have used in the last equation that $\sum _{j=1}^{k} \bar X_j - p_j = 0$

The latter equation can be rewritten in terms of matrix notation as
\begin{eqnarray}
\chi^2 & = & n (\bar X^* - p^* )^T (\Sigma^*) ^{-1}  (\bar X^* - p^* )
\end{eqnarray}

The Central limit theorem states that $Y_n= \sqrt n (\Sigma^*) ^{-1/2} (\bar X^* - p^* )$ converges in distribution to a normal variable $N_{k-1}(0,I_{k-1})$. The $\chi^2$ statistic given by $
(Y_n)^T Y_n$ converges in distribution to $\chi^2_{\infty} = N_{k-1}(0,I_{k-1})^T N_{k-1}(0,I_{k-1})$. $\chi^2_{\infty}$  is the sum of the squares of $k-1$ independent standard normal random variables, which is a chi square distribution with $k-1$ degree of freedom. This concludes the first proof. \qed

\subsection{Second Proof: Cochran theorem}
The second proof relies on the Cochran theorem. The start is the same. Since $\mathbb{E}[X_i] = p$, the central limit theorem implies 
\begin{equation}\label{TCL}
\sqrt n ( \bar{X}_n - p ) \xrightarrow[n \rightarrow \infty]{\text{d}} N_k(0, \Sigma)
\end{equation}

where $N_k$ denotes the $k$ multi dimensional normal as always in this paper. Let us denote by $Z$ the standard reduced Gaussian variable corresponding to the basis implied by the $k$ multi dimensional normal $N_k(0,I_k)$. 
 
If we apply the Cochran theorem with the projection on the sub vectorial space $F$ spanned by $\sqrt{p} = (  \sqrt{p_1}, ..., \sqrt{p_k})^T$ (whose norm is obviously 1), we have that the projection matrix on $F$ is given by 
$$P_F = \sqrt{p}(\sqrt{p}^T\sqrt{p})^{-1} \sqrt{p}^T = \sqrt{p} \sqrt{p}^T$$  and that the projection on the orthogonal of $F$ denoted by $F^{\bot}$ is given by $P_{F^{\bot} } = I - P_F = I - \sqrt{p} \sqrt{p}^T$.

Since $F$ is spanned by one single vector, its dimension is 1, while its orthogonal, $F^{\bot}$, is of dimension $k-1$.
The Cochran theorem states that the projection on $F^{\bot}$ of $Z$ follows a normal whose distribution is given by $N(0,P_{F^{\bot}})$, and that the squared norm of the projection on $F^{\bot}$ follows a chi square distribution of dimension $k-1$. 

We can notice that if we define $ \Gamma  = \text{diag} (p)$ and $A_n = \sqrt n \Gamma^{-1/2} (\bar X - p )$. 
The Chi-squared statistics can be rewritten as the norm of the stochastic vector $A^n$ since
\begin{eqnarray}
\chi^2 & = (A_n)^T A_n
\end{eqnarray}

Using equation (\ref{TCL}), we also know that $A_n$ converges in distribution to $N(0,  \Gamma^{-1/2}  \Sigma \Gamma^{-1/2})$. Since $\Sigma = \Gamma - p p ^T$, we have that 
\begin{eqnarray}
\Gamma^{-1/2}  \Sigma \Gamma^{-1/2} &= &  I_k -\Gamma^{-1/2}  (p p ^T) \Gamma^{-1/2}  = I_k - (\Gamma^{-1/2}  p)  (\Gamma^{-1/2} p)^T  \\
&=&  I_k -\sqrt{p} \sqrt{p}^T = P_{F^{\bot}}
\end{eqnarray}

This states that $A_n$ converges in distribution to the projection of Z on $F^{\bot}$. The statistics, $\chi^2$, that is the squared norm of $A_n$ converges in distribution to the squared norm of the projection of Z on $F^{\bot}$, whose distribution is a chi square distribution of dimension $k-1$. This proves that the statistics, $\chi^2$, converges in distribution to a chi square distribution of dimension $k-1$. \qed

\subsection{Third Proof: Sylvester theorem and eigen values}
The third proof relies on Sylvester theorem and is based on an explicit computation of the eigen values of the associated covariance matrix. If we write $Z$ the $k$ vector with coordinates given by $Z_i = \frac{N_i-np_i}{\sqrt{np_i}} = \sqrt n \frac{ N_i / n - p_i}{ \sqrt p_i} $, the central limit theorem states that vector $Z$ converges in distribtuion to $N(0,\Omega)$, a multivariate normal distribution, whose covariance matrix is given by

$$\Omega=\text{Cov}(Z)=\begin{pmatrix}
1-p_1 & -\sqrt{p_1 p_2} & \cdots \\
-\sqrt{p_1 p_2} & 1-p_2 & \cdots \\
\vdots & \vdots & \ddots
\end{pmatrix}.$$

We can compute explicitly the characteristic polynomial of this matrix. The Sylvester theorem states that $\text{Det}(I_{k} - c r) = 1 - r^Tc$ for any $c, r \in \mathbb{R}^k$ ($k$ dimensional vectors). Therefore, a direct application of Sylvester theorem shows that $\text{Det}(\Omega-\lambda I)=(1-\lambda)^{k-1}\lambda$  as $\Omega=I-pp^T$ for $p=(\sqrt{p_1},\sqrt{p_2},\dots)$ and $\text{Det}(\Omega-\lambda I)=(1-\lambda)^{k} \text{Det}(I_{n} - \frac{1}{1-\lambda} pp^T )$ . This implies that $\Omega$ has $k-1$ eigenvalues that are 1 and one that is 0 and that the distribution is really $k-1$ dimensional embedded in $k$ dimensions. In particular there is a rotation matrix $A$ that makes
 
$$A \Omega A^T=\begin{pmatrix}
1 		& 0 		& \ldots  	& 0 		& 0 		\\
0 		& 1 		& \ldots  	& 0 		& 0 		\\
\vdots 	& \vdots & \ddots  & \vdots 	& \vdots 	\\
0 		& 0 		& 0		& 1 		& 0 		\\
0 		& 0 		& 0		& 0 		& 0 		
\end{pmatrix} 
= \begin{pmatrix} 
I_{k-1} & 0_{k-1,1} \\
0_{1,k-1}	    & 0_{1,1}
\end{pmatrix} 
$$
where $0_{n,m}$ is the $n$ rows, $m$ columns matrix filled with 0. 

Denote $W = AZ \sim N_k(0,A\Omega A^T)$. Then $W$ is a vector $(W_1, W_2, \dots, W_{k-1}, 0 )$ of iid. $\mathcal N(0,1)$ Gaussians with only $k-1$ non null coordinates (the first $k-1$ coordinates). The function $f(Z) = Z_1^2 + Z_2^2 + \dots$ is the norm $\|Z\|_2^2$, and hence it is invariant if we rotate its argument. This means $f(Z) = f(AZ) = f(W) = W_1^2 + W_2^2 + \dots + W_{k-1}^2$ is Chi-square distributed with $k-1$ degree of freedom\qed

\subsection{Fourth Proof: Characteristic function}
Let us use characteristic function. The $\chi^2$ statistic' s characteristic function is:
\begin{eqnarray}
\phi_\chi^2( t ) &=& \mathbb E \left[  e^{\mathrm i t \chi^2 }\right] = \mathbb E \left[  e^{\mathrm i t  \sum _{j=1}^{k} \frac{ (N_j - n p_j)^2 }{ n p_j}  }\right] 
\end{eqnarray}

Since $ \sum_{ j=1 \ldots k} N_j - n p_j =0$, we can reduce the sum to $k-1$ terms and show the real quadratic form as follows: 

\begin{eqnarray}
\phi_\chi^2( t ) &=& \mathbb E \left[  e^{  \mathrm{i} t  \left( \sum _{j=1}^{k-1} \frac{ (N_j - n p_j)^2 }{ n p_j} + \frac{ ( \sum _{j=1}^{k-1}  N_j - n p_j )^2 }{  n p_k}  \right) }\right] \\
&=& \mathbb E \left[  e^{  \frac{ \mathrm{i} t } {n}  \left(  (N^{*} - n p^{*})^T  (\Sigma^{*})^{-1}  (N^{*} - n p^{*})   \right) }\right] 
\end{eqnarray}

where  $N^{*} = (N_1, \ldots, N_{k-1})^T$ is a $k-1$ stochastic vector and $(\Sigma^*) ^{-1}$ the $k-1$ squared symmetric matrix is given by
\begin{equation}
(\Sigma^*) ^{-1} = \left(
\begin{array}{l l l l}
{1/p_1}+{1/p_k} 	& {1/p_k}			& \ldots 		&  {1/p_k}		\\ 
{1/p_k}			& {1/p_2}+{1/p_k} 	& \ldots 		& {1/p_k}		\\ 
\vdots 			& \vdots 				& \ddots 		& \vdots 			\\
{1/p_k}			& {1/p_k}			& \ldots 		& {1/p_{k-1}}+{1/p_k} 	
\end{array} \right)
\end{equation}

The variance matrix of $N^*$ is easily computed and given by :
\begin{equation}
\Sigma^*  = \left(
\begin{array}{l l l l}
{p_1 (1-p_1)}	& {-p_1 p_2}	& \ldots 		& {-p_1 p_{k-1}}		\\
{-p_1 p_2}		& {p_2 (1-p_2)}	& \ldots 		& {-p_2 p_{k-1}}		\\
\vdots 			& \vdots 		& \ddots 	& \vdots 		  	\\
{- p_1 p_{k-1} } 		& {- p_2 p_{k-1} }  	& \ldots 		& {p_{k-1} (1-p_{k-1})}
\end{array} \right)
\end{equation}

It is remarkable that looking at the characteristic function provides the inverse of $\Sigma^*$ without effort as one can validate that $(\Sigma^*) ^{-1}  \Sigma^* = I_{k}$.
The Central limit theorem for Binomial (which is also referred to as De Moivre-Laplace's theorem) states that $N^{*} - n p^{*}  \xrightarrow[n \rightarrow \infty]{\text{d}} N( 0, \text{Var}(N^{*}))= N( 0, \Sigma^* )$.  Hence, $ (\Sigma^{*})^{-1/2}  (N^{*} - n p^{*})  \xrightarrow[n \rightarrow \infty]{\text{d}} N( 0, I_{k-1})$. Let us denote by $Z = (Z_1, \ldots, Z_{k-1})^T$ the corresponding $k-1$ dimensional standard normal distribution. Taking the limit in the characteristic function thanks to Lévy's continuity theorem, we have therefore that

\begin{eqnarray}
\phi_\chi^2( t )  &\xrightarrow[n \rightarrow \infty] {\text{d}}  & \mathbb E \left[   \prod  _{ j=1, \ldots, k-1} e^{\mathrm i t  Z_j^2} \right]  = \left(  \mathbb E \left[  e^{\mathrm i t  U^2 } \right]\right)^{k-1}
\end{eqnarray}

where $U$ is a standard normal $N(0,1)$. 
We have $\mathbb E \left[  e^{\mathrm i t  U^2 } \right] = \frac{1}{1-2 i t}$ which leads to $\phi_X( \chi^2 ) \underset{n \rightarrow \infty}{\rightarrow} \frac{1}{(1-2 it ) ^{k-1}}$ which concludes the proof as the characteristic function of a $\chi^2(k-1)$ is precisely $ \frac{1}{(1-2 it ) ^{k-1}}$ \qed

\subsection{Fifth Proof: Projection matrix and Pythagoras theorem}
Even if this proof is relying on similar argument as the Cochran theorem, it slightly differs and have been first shown in \cite{Hunter_2015}. We define $\Gamma = diag(p)$. The central limit theorem states that 
\begin{equation}
\sqrt n \Gamma^{-1/2} ( \bar X - p ) \xrightarrow[n \rightarrow \infty]{\text{d}} N_{k}(0, \Gamma^{-1/2} \Sigma \Gamma^{-1/2})  
\end{equation}

Noticing that $\Sigma = \Gamma - p p ^T$, we  have 
\begin{equation}
\Gamma^{-1/2} \Sigma \Gamma^{-1/2} = I_k - \Gamma^{-1/2} p p^T \Gamma^{-1/2} = I_k- \sqrt p \sqrt{ p} ^T
\end{equation}

Using the linearity and the commutativity property of the trace, we have 
\begin{equation}
\text{Trace}( \Gamma^{-1/2} \Sigma \Gamma^{-1/2} ) = \text{Trace}(I) - \text{Trace}( \sqrt{ p}^T  \sqrt p ) = k-1
\end{equation}
since $ \sqrt{ p}^T  \sqrt p = 1$. We can also notice that 
\begin{equation}
( I_k- \sqrt p \sqrt{ p} ^T)^2 =  I_k - 2  \sqrt p \sqrt{ p} ^T +   \sqrt p \sqrt{ p} ^T =  I_k- \sqrt p \sqrt{ p} ^T
\end{equation}
which means that $ \Gamma^{-1/2} \Sigma \Gamma^{-1/2} $ is an idempotent matrix. Denoting by $A_n = \sqrt n \Gamma^{-1/2} ( \bar X - p )$, we can notice that $\chi^2$ is the squared norm of $A_n$: $\chi^2= A_n^T A_n = \norm{A_n}^2$. Since $ I_k- \sqrt p \sqrt{ p} ^T$ is idempotent, it is a projection matrix of rank equal to its trace: $k-1$. We can conclude using the following lemma found in \cite{Hunter_2015}. 

\begin{lemma}
Suppose $P$ is a projection matrix. Then if $Z \sim N_k (0, P )$, $Z^T Z \sim \chi^2(r)$ where r is the trace of $P$ or equivalently the number of eigen values of $P$ equal to 1 or equivalently the number of eigen values of $P$  non equal to 0.
\end{lemma}

\qed.

\subsection{Sixth Proof: Generic induction with De Moivre-Laplace theorem}
The sixth proof differs from previous ones in the spirit as it proves it using generic induction method. It also uses a weaker form of the central limit theorem (the De Moivre-Laplace theorem) that provides the asymptotic distribution for a binomial distribution. The goal here is to prove the following proposition.

\begin{proposition}
if for $k>2$, the $Q_k$ statistics given by $Q_k= \sum _{j=1}^{k} \frac{ (N_{j,k} - n p_{j,k})^2 }{ n p_{j,k}}$
 converges in distribution to a $\chi^2(k-1)$ with  $\sum _{j=1}^{k} p_{j,k}=1$ and $\sum _{j=1}^{k} N_{j,k}=n$, \newline
\indent then the $k+1$ statistic given by 
$Q_{k+1}= \sum _{j=1}^{k+1} \frac{ (N_{j,k+1} - n p_{j,k+1})^2 }{ n p_{j,k+1}}$
converges in distribution to a $\chi^2(k)$ with $\sum _{j=1}^{k} p_{j,k+1}=1$ and $\sum _{j=1}^{k} N_{j,k+1}=n$ 
\end{proposition}

To emphasize the fact that the probabilities and partition of n at rank $k$ are different from the one at rank $k+1$, we have use the underscore notation for $p_{j,k}$ and $N_{j,k}$. However, when the underscore notation will be obvious, we will drop it to make the computation more readable. We will anyway warn the reader to make sure we do not loose him or her.

\begin{proof}
For k=2,  $Q_k$ writes as
\begin{eqnarray}
Q_{k} =  \frac{ (N_{1,2} - n p_{1,2})^2 }{ n p_{1,2}} + \frac{ (N_{2,2} - n p_{2,2})^2 }{ n p_{2,2}}  =  \left[ \frac{ N_{1,2} - n p_{1,2}}{ \sqrt{ n p_{1,2} (1- p_{1,2}) } } \right]^2
\end{eqnarray}

since $ p_{1,2}+  p_{2,2}=1$ and $ N_{1,2} + N_{2,2}=n$.

De Moivre-Laplace's theorem states that $X_n= \frac{ N_{1,2} - n p_{1,2}}{ \sqrt{ n p_{1,2} (1- p_{1,2}) } } $ 
converges in distribution to a standard normal distribution $N(0,1)$ since $N_{1,2} \sim Bin(n, p_{1,2})$. Hence, $Q_{k} \xrightarrow[n \rightarrow \infty]{\text{d}} \chi^2(1)$.

Suppose the property to prove is true for $n=k$. We can create a general induction between $Q_{k}$ and $Q_{k+1}$ as follows:
\begin{eqnarray*}
Q_{k+1} &= & \sum _{j=1}^{k+1} \frac{ (N_{j,k+1} - n p_{j,k+1})^2 }{ n p_{j,k+1}} \\
&= &  \sum _{j=1}^{k-1} \frac{ (N_{j,k+1} - n p_{j,k+1})^2 }{ n p_{j,k+1}} 
+ \frac{  (N_{k,k+1}+ N_{k+1,k+1}- n (p_{k,k+1}+p_{k+1,k+1}))^2 }{ n (p_{k,k+1}+p_{k+1,k+1})} + U^2 \\
& = &  \sum _{j=1}^{k} \frac{ (N^{'}_{j,k} - n p^{'}_{j,k})^2 }{ n p^{'}_{j,k}}  + U^2 \\
& = &  Q^{'}_{k} + U^2 \\
\end{eqnarray*}	

\vspace{-0.5cm}
where we have used the following notations
\begin{eqnarray*}
N^{'} &=& ( N_{1,k+1}, \ldots, N_{k-1,k+1}, N_{k,k+1}+N_{k+1,k+1} ), \\
p^{'} &=& (p_{1,k+1}, \ldots, p_{k-1,k+1}, p_{k,k+1}+p_{k+1,k+1}), \\
U^2 &=& \frac{ (N_{k,k+1} - n p_{k,k+1})^2 }{ n p_{k,k+1}} + \frac{ (N_{k+1,k+1} - n p_{k+1,k+1})^2 }{ n p_{k+1,k+1}}   - \frac{  (N_{k,k+1}+ N_{k+1,k+1}- n (p_{k,k+1}+p_{k+1,k+1}))^2 }{ n (p_{k,k+1}+p_{k+1,k+1})},\\
Q^{'}_{k} &=& \sum _{j=1}^{k} \frac{ (N^{'}_{j,k} - n p^{'}_{j,k})^2 }{ n p^{'}_{j,k}} 
\end{eqnarray*}	

By assumption, $Q^{'}_{k}  \xrightarrow[n \rightarrow \infty]{\text{d}} \chi^2(k-1)$. 
The last remaining part is to prove that $Q^{'}_{k}$ and $U$ are independent and that $U^2 \sim \chi^2(1)$. 
In the following to make notation lighter, we will drop the lower index $._{k+1}$. 
Let us denote by $T_k =N_{k} - n p_{k} $, $T_{k+1} = N_{k+1} - n p_{k+1}$, $q_k=\sqrt{ p_{k}}$, $q_{k+1} = \sqrt{p_{k+1}}$, a straight computation leads to
\begin{eqnarray*}
U^2 &=& \frac{1}{ n (p_{k}+p_{k+1})} \left[ \frac{p_{k+1} T_k^2 }{ p_{k} } + \frac{ p_{k}  T_{k+1}^2 }{  p_{k+1}  }  - 2   T_k T_{k+1} \right]  \\
& =& \left[ \frac{1}{ \sqrt{ n (p_{k}+p_{k+1})}}  \right]  ^2   \left[ \frac{ q_{k+1}  T_{k}  }{q_{k}} -\frac{ q_{k}  T_{k+1}  }{q_{k+1}}  \right]  ^2 \\
& =& \left[ \frac{ p_{k+1} N_{k} - p_{k} N_{k+1} }{ \sqrt{ n p_{k} p_{k+1} (p_{k}+p_{k+1})}}  \right]  ^2   \\
& = & V^2
\end{eqnarray*}	

where $V=\frac{ p_{k+1} N_{k} - p_{k} N_{k+1} }{ \sqrt{ n p_{k} p_{k+1} (p_{k}+p_{k+1})}} $. 
De Moivre-Laplace's theorem states that $N = (N_1,...,N_{k+1})$ converges in distribution to a Gaussian vector, or consequently that $V$ converges in distribution to a normal distribution, 
as $V$ is a linear combination of the coordinate of $N$.  $V$'s mean is simple to calculate and equal to 0 since \useshortskip 
\begin{equation*}
\mathbb{E}[  p_{k+1} N_{k} - p_{k} N_{k+1}  ] = n ( p_{k+1} p_{k}- p_{k} p_{k+1}).
\end{equation*}
\nr[-1.2cm] $V$'s variance is simple to calculate and equal to 1 since \useshortskip 
\begin{equation*}
\text{Var}[  p_{k+1} N_{k} - p_{k} N_{k+1}  ] = n \left[ p_{k+1}^2 (p_{k} (1-p_k) )+  p_{k}^2 (p_{k+1} (1-p_{k+1})) - 2  p_{k}^2 p_{k+1}^2  \right].
\end{equation*}
\nr[-1.2cm] Hence $V  \xrightarrow[n \rightarrow \infty]{\text{d}} N(0,1)$. And $U^2$ converges in distribution to a $\chi^2(1)$ distribution. \\

\vspace{-0.2cm} The final part is to prove the independence of $U$ and $Q^{'}_{k}$ or equivalently the independence of $L= p_{k+1} N_{k} - p_{k} N_{k+1} $ and $Q^{'}_{k}$. $Q^{'}_{k}$ is composed of coordinates of $N^{'}$. So it is is sufficient to prove that $L$ is independent of all the coordinates of $N^{'}$. Both $N^{'}$ and $L$ are coordinates of a Gaussian vector. So, their independence is equivalent to a null covariance. Let us compute. For any $j \leq k+1$, $\text{Cov}(N_j, L) =  n (  p_{k+1} p_{j} p_{k} - p_{k} p_{j} p_{k+1} ) = 0$. The covariance with the last coordinate of $N^{'}$ is also null as it is $\text{Cov}(N_{k}+N_{k+1}, L ) = \text{Cov}(N_{k},L) +  \text{Cov}(N_{k+1},L) = 0 + 0 = 0$. This concludes the proof. \qed
\end{proof}

\subsection{Seventh Proof: Connection with Poisson variables and geometry}
This proof is due to \cite{Fisher_1922} and relives on geometry arguments. In the sequel, we shall write $k$ to be an integer, $I_k$ the identity matrix of order $k$
and $Z = (Z_1, Z_2, \ldots , Z_k)$ a random vector with a multi dimensional standard normal distribution
$N_k(0, I_k)$ of order dimension $k$. We shall assume the components of $Z$ to be independent. 

Fisher used a geometric argument to determine the distribution of the
random variable
$$
U := Z^2_1 + Z^2_2 + \ldots + Z^2_k
$$

The values $ (Z_1, Z_2, \ldots , Z_k)$  of any given sample of $Z$ can be interpreted as 
co-ordinates of a point $P$ in the $k$-dimensional Euclidean space $\mathbb{R}^k$. 

The Euclidean distance of $P$ from the origin $O$  is written as $U$ and it represents its $L^2$ norm defined by $U = \norm{ OP}^2 $. One property of the Euclidean distance is to be unchanged by any rotation of the co-ordinates orthonormal axes. 

The joint probability density function of the components of $Z$ should therefore be proportional to $e^{-\norm{OP}^2/2}$ and should remain constant on the $k$-dimensional hypersphere with radius $\sqrt u = \norm{OP} $.
A consequence is that the density $f_U(u)$ shall be obtained by integrating it between the two hyperspheres
with radius $u$ and $u + du$

\begin{eqnarray}\label{FisherAssumption}
f_U(u) &= & c e^{-\norm{OP}^2/2} \frac{d}{d\norm{OP}} \norm{OP}^k  \\
&= & c e^{u/2} \frac{d}{du} u^{k / 2}
\end{eqnarray}

We shall also impose that there is $c$ a suitable normalization constant to have a density summed to 1. 
Hence, the constant shall be the following:
\begin{eqnarray}
f_U(u) &= & \frac{1}{2^{r/2} \Gamma(k /2)}  e^{-u/2} u^{k /2-1}
\end{eqnarray}
Hence, $U \sim \chi^2(k)$.

In the particular case of $s \leq  k$ linear (independent) constraints for the components of $Z$, we shall be able to generalize previous results.
Each constraint defines a hyper-plane of $\mathbb{R}^k$ containing $O$, say
$\pi$, and $Z$ is forced to belong to it. The intersection of the generic
hyper-sphere of $\mathbb{R}^k$ with $\pi$ is a hyper-sphere of the space   $\mathbb{R}^{k-1}$. The result
of the $s$ linear constraints will create a hypersphere of $\mathbb{R}^{k-s}$ from a
generic hypersphere of $\mathbb{R}^{k}$. We can apply previous reasoning with the adaptation or replacing $k$ by $k-s$. 
In other words, if there are $s \leq k$ independent linear constraints for $Z$, one has
\begin{equation}
U = \sum_{i=1}^{r} Z_i^2  \sim \chi^2(k-s).
\end{equation}

Let us consider $V = (V_1, V_2, \ldots, V_k)$ consisting of $k$ independent random
variables with Poisson distribution with intensity parameter $n p_{0,i}$  for $i = 1, 2\ldots,  k$ : 
$V_i \sim \Pi(n p_{0,i})$. 

With $(n_1, n_2, \ldots, n_k, n ) \in \mathbb{N}^{k+1}$ and $N = \sum_{i=1}^k V_i$, at first, one has:
\begin{eqnarray}\label{poisson1}
P  (V_1=n_1, V_2=n_2, \ldots, V_k=n_k, N=n) &=& \prod_{i=1}{k}\frac{(n p_{0,i}) ^{n_i}}{n_i ! }e^{-n p_{0,i}} \\
& = & e^{-n} n ^n  \prod_{i=1}{k}\frac{( p_{0,i}) ^{n_i}}{n_i ! }
\end{eqnarray}

Note that, although in the left-hand side of (\ref{poisson1}), $N$ appears as a random variable,
the joint distribution is singular because $N =  \sum_{i=1}{k} N_i$. 
The marginal distribution of $N$ is easy to determine because $N$ is the sum of $k$
independent random variables with Poisson distribution and with mean $\mathbb{E}(N) = n$. Hence
$N  \sim \text{Poisson}(n)$. From this fact and from equation (\ref{poisson1}), one gets the conditional probability:

\begin{eqnarray}
P  (V_1=n_1, V_2=n_2, \ldots, V_k=n_k \vert N=n) &=& \frac{ n ! } {n_1 ! n_2! \ldots n_k!} \prod_{i=1}{k}(p_{0,i}) ^{n_i}
\end{eqnarray}

This proves that the distribution of $V$ , under the condition $N = n$ is the same as the one defined by $N = (N_1, N_2, \ldots, N_k)$, the $k$ binomial variables corresponding to the $k$ categorical variables. 
Let us write
$$
\tilde{Z}_i = \frac{V_i - n p_{0,i}}{\sqrt{n p_{0,i}} }
$$
The Central Limit theorem (for Poisson variable with growing intensity parameter)  states that the renormalized Poisson variables $\tilde{Z}_i$ converges to standard normals:
$$
\tilde{Z}_i = \frac{V_i - n p_{0,i}}   {\sqrt{n p_{0,i}}}  \xrightarrow[n \rightarrow \infty]{\text{d}} N(0, 1)
$$
These Gaussian variables are not independent as the variables $\tilde{Z}_i$ have to satisfy the constraint
$$
\sum_{i=1}^{k} \sqrt{n p_{0,i}} \tilde{Z}_i = 0
$$

The above reasoning states that the squared norm of these Gaussian variables should have a $\chi(k-1)$ distribution as they have to satisfy one linear constraints. Wrapping everything up , we have that 
$$
\sum_{i=1}^{k} \tilde{Z}_i^2 = \sum_{i=1}^{k}  \frac{ (N-i - {n p_{0,i} ) }} {n p_{0,i}}   \xrightarrow[n \rightarrow \infty]{\text{d}} \chi^2 (k-1)
$$
which concludes the proof.

\qed

\section{Graphical interpretation of the test}
In the different proofs, we have seen that  asymptotically, the variables $\frac{N_j - n p_j}{\sqrt{ n p_j (1-pj)}}$ are normally distributed. This has given us the idea, for a two by two (2x2) contingency table, to not only draw the table with color and size to give an intuition of the relationship between the two variables but also to provide a confidence interval to illustrate whether the second categorical variable is close or not to the first one in the sense of the Pearson correlation test.

\begin{figure}[h]
\begin{center}
\includegraphics[width=11.5cm]{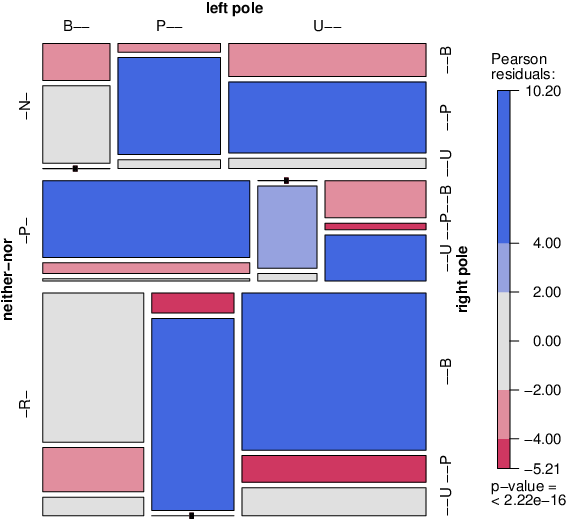} 
\caption{Standard mosaic for a two-way contingency table. The color and the size provides information about the values of the contingency table. On the right side, the Pearson scale provides information about the residuals. However no confidence interval is provided.}
\label{figure1}
\end{center}
\end{figure}

Indeed, in statistical graphics, mosaic display, attributed to \cite{Hartigan_Kleiner_1981}, is a graphical method to show the values (cell frequencies) in a contingency table cross-classified by one or more factors. Figure \ref{figure1} shows the basic form of a mosaic display for a two-way table of individuals. Canonical examples can be found in 
\cite{Friendly_1992}, \cite{Friendly_1994} and \cite{Friendly_2002}. Mosaic displays have become the primary graphical tool for visualization and analysis categorical data in the form of n-way contingency tables. Although they provide Pearson residuals, it appeared to us that adding a confidence interval could increased readability. This is what we have achieved as follows. The full source code used can be found in \href{https://github.com/ericbenhamou/Pearson-mosaicplot}{github}.

\begin{figure}[h!]
\begin{center}
\includegraphics[width=8cm]{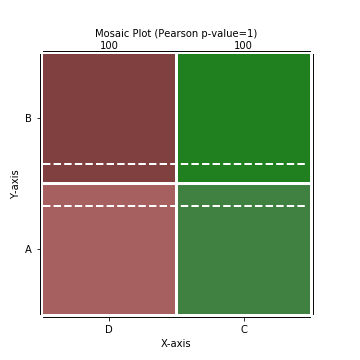} 
\caption{The cases where all values are equal to 50. Each individual square is of same size. The dark red and dark green square lower side ore limits are within the confidence interval. They lie precisely at the middle of the confidence interval indicating that the two categorical variables are very very similar or from the same distribution (with regards to the Pearson independence test). We can conclude that these two categorical variables are statistically not independent. This case is our benchmark.}
\label{mosaic-initial}
\end{center}
\end{figure}

The method works as follows:
We first compute the total values of the our contingency table. This is the sum of all the values inside the table. We call this N. We estimate the probability $\hat{p}$ over all categories as follows:  $\hat{p}$ = values of interest / $N$. Note that this probability is estimated on all data (under the assumption that both categories are from the same distribution). A confidence interval for the probability estimator is easily calculated from the standard normal distribution:  $\delta_{p} = \text{pdf( quantile )} \sqrt( \hat{p} (1 - \hat{p} ) \ N$ where pdf stands for the traditional probability density function. We then just to add this confidence interval to our table centered around our estimated probability. In figure \ref{mosaic-initial}, we provide the mosaic plot for an initial table given by ideal values provided in \ref{initial-table}.

\begin{table}[h!]
\begin{center}
\begin{tabular}{| l | l | l |}
\hline
 	& C 		& D  	\\ \hline
A 	& 50 	& 50		\\ \hline
B 	& 50 	& 50		\\ \hline
\end{tabular}
\caption{initial data for the contingency table}
\label{initial-table}
\end{center}
\end{table}

\begin{figure}[h!]
\begin{center}
\includegraphics[width=4.5cm]{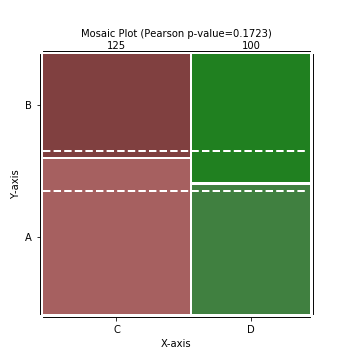} 	\includegraphics[width=4.5cm]{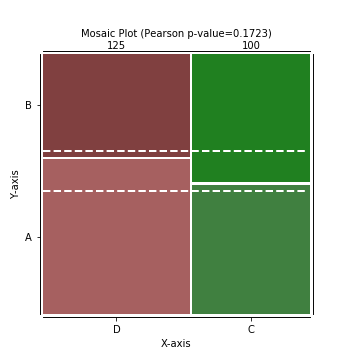}  \\
\includegraphics[width=4.5cm]{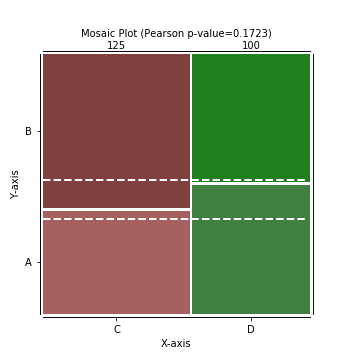} 	\includegraphics[width=4.5cm]{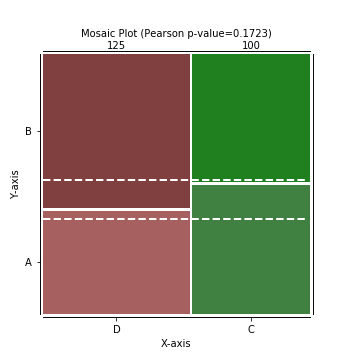} 
\caption{We plot the four cases where we change one value to 75 and keep the other ones to 50. We start departing from the same distribution. However, in all four cases, the modified value is still within the confidence interval. It is further away from the middle of the large square but still within the confidence interval. It is worth noting that the estimated probability is computed on the full data. Hence, the confidence interval is shifted either upward or downward depending on the case. As we always plot in the first column the categorical variable with more observation, the graphics are the same on the right and left}
\label{mosaic-75}
\end{center}
\end{figure}

\begin{figure}[h!]
\begin{center}
\includegraphics[width=4.5cm]{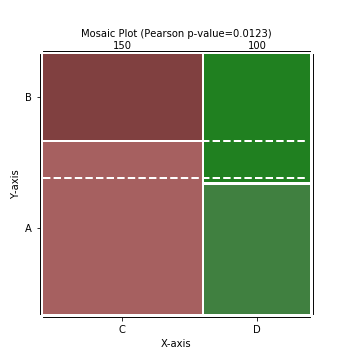} 	\includegraphics[width=4.5cm]{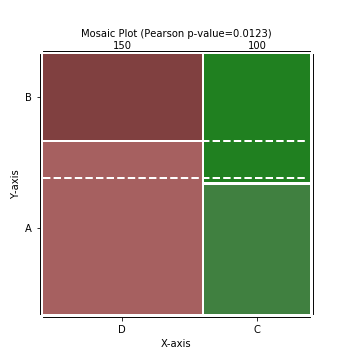}  \\
\includegraphics[width=4.5cm]{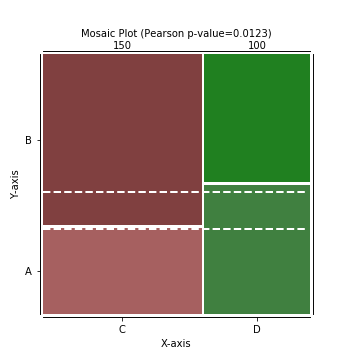} 	\includegraphics[width=4.5cm]{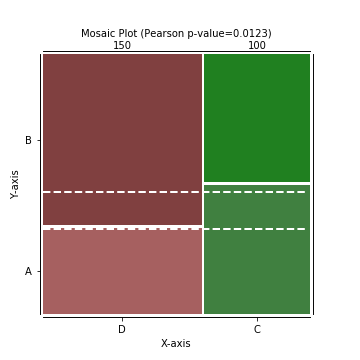} 
\caption{The four cases where we have change one value to 100 and keep the other to 50. In this case, we are slightly outside from the interval confidence interval. This indicates in particular that if we had taken a higher error of type I, the test could have been successful.}
\label{mosaic-100}
\end{center}
\end{figure}

\begin{figure}[h]
\begin{center}
\includegraphics[width=4.5cm]{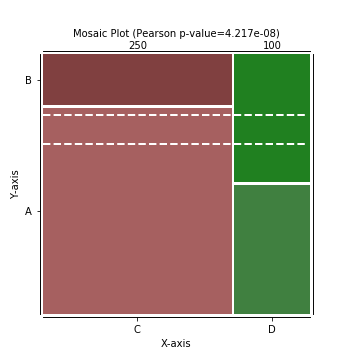} 	\includegraphics[width=4.5cm]{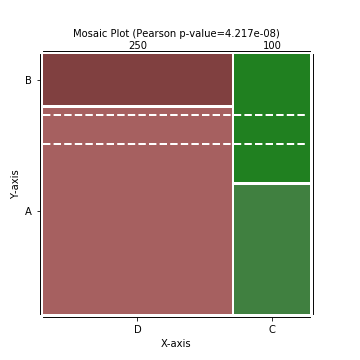}  \\
\includegraphics[width=4.5cm]{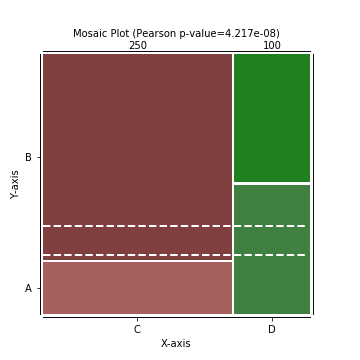} 	\includegraphics[width=4.5cm]{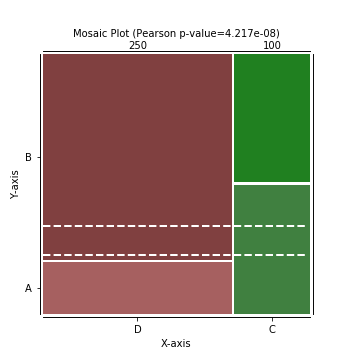} 
\caption{The four cases where we have change one value to 200 and keep the other to 50. We should note that these figures are very different from figure \ref{mosaic-100}. In our case, we are much more outside from the confidence interval  than in the case of figure \ref{mosaic-100}. This is precisely the interest of this graphical representation to be able to compare and dissociate the two cases. Values for these graphics are provided in table\ref{table-200}.}
\label{mosaic-200}
\end{center}
\end{figure}

We then change each of the four values with the following new values: 75, 100 and 200 to visualize the impact of a variable that becomes more and more independent. Corresponding graphics are provided in \ref{mosaic-75}, \ref{mosaic-100} and \ref{mosaic-200} with corresponding tables \ref{table-75}, \ref{table-100} and \ref{table-200}. In the initial case (provided in figure {mosaic-initial} with corresponding values given in table \ref{initial-table}), each individual square is of equal size and equal to one fourth of the total square. This is logical as all values are equal to 50. Since the dark red and dark green square lower side are within the confidence interval, at the middle of the confidence interval, we can graphically visualize that the two categorical are from the same distribution (with regards to the Pearson independence test). We will use this case as a benchmark and progressively change the second categorical variables to make it independent (or in the sense of the Pearson test, the second categorical variable will not come from the same distribution from a statistical point of view). The case with a modified value of 75 provides some hindsight about the graphical interpretation of the Pearson test). It is given by table \ref{mosaic-75}. First of all the lower or upper side of the concerned squares lie within the confidence interval. This means that from a statistical point of view, the two categorical variables are from the same distribution. As we see that these two categorical variables are drifting away one from another, we can see that compared to our benchmark the two categorical variables are slightly different, though from the same distribution! Secondly, we see that the confidence interval is shifted upward or downward depending on the case. Because we always first show the categorical variable with more observation for all our figures (\ref{mosaic-75}, \ref{mosaic-100} and \ref{mosaic-200}), the left and right column are the same, though x axis indexes are permuted. 
Figure \ref{mosaic-75} is an example of Pearson test where we fail to reject hypothesis H0 meaning we fail to reject H0: the two categorical variables are from the same distribution. Figure \ref{mosaic-100} is informative. In this case, we are slightly outside the confidence interval indicating that in this particular case, we can reject H0: 
the two categorical variables are from the same distribution.  However, this is way different from the case 
of figure \ref{mosaic-200} where we are much more outside the confidence interval. 
Figure \ref{mosaic-100} corresponds to a modified value of 100 versus figure \ref{mosaic-200}  that corresponds to a modified value of 200. 
This is precisely the interest of this graphical representation. We can measure by how much we are from failing to reject hypothesis H0.  In all tables that provide the modified value (tables \ref{table-75}, \ref{table-100} and \ref{table-200}), we have emphasized the modified value by writing it in bold.

\begin{table}[h!]
\begin{center}
\begin{tabular}{l l l l l l}
                      &     C                & D                    &                       & C                    & D                     \\ \cline{2-3} \cline{5-6} 
\multicolumn{1}{l|}{A} & \multicolumn{1}{l|}{\textbf{75}} & \multicolumn{1}{l|}{50} & \multicolumn{1}{l|}{} & \multicolumn{1}{l|}{50} & \multicolumn{1}{l|}{\textbf{75}} \\ \cline{2-3} \cline{5-6} 
\multicolumn{1}{l|}{B} & \multicolumn{1}{l|}{50} & \multicolumn{1}{l|}{50} & \multicolumn{1}{l|}{} & \multicolumn{1}{l|}{50} & \multicolumn{1}{l|}{50} \\ \cline{2-3} \cline{5-6} 
                      &                       &                       &                       &                       &                       \\
                      &     C                & D                    &                       & C                    & D                     \\ \cline{2-3} \cline{5-6} 
\multicolumn{1}{l|}{A} & \multicolumn{1}{l|}{50} & \multicolumn{1}{l|}{50} & \multicolumn{1}{l|}{} & \multicolumn{1}{l|}{50} & \multicolumn{1}{l|}{50} \\ \cline{2-3} \cline{5-6} 
\multicolumn{1}{l|}{B} & \multicolumn{1}{l|}{\textbf{75}} & \multicolumn{1}{l|}{50} & \multicolumn{1}{l|}{} & \multicolumn{1}{l|}{50} & \multicolumn{1}{l|}{\textbf{75}} \\ \cline{2-3} \cline{5-6} 
\end{tabular}
\caption{Data for the various figures \ref{mosaic-75}. We have changed one of the value to 75}
\label{table-75}
\end{center}
\end{table}

\begin{table}[h!]
\begin{center}
\begin{tabular}{l l l l l l}
                      &     C                & D                    &                       & C                    & D                     \\ \cline{2-3} \cline{5-6} 
\multicolumn{1}{l|}{A} & \multicolumn{1}{l|}{\textbf{100}} & \multicolumn{1}{l|}{50} & \multicolumn{1}{l|}{} & \multicolumn{1}{l|}{50} & \multicolumn{1}{l|}{\textbf{100}} \\ \cline{2-3} \cline{5-6} 
\multicolumn{1}{l|}{B} & \multicolumn{1}{l|}{50} & \multicolumn{1}{l|}{50} & \multicolumn{1}{l|}{} & \multicolumn{1}{l|}{50} & \multicolumn{1}{l|}{50} \\ \cline{2-3} \cline{5-6} 
                      &                       &                       &                       &                       &                       \\
                      &     C                & D                    &                       & C                    & D                     \\ \cline{2-3} \cline{5-6} 
\multicolumn{1}{l|}{A} & \multicolumn{1}{l|}{50} & \multicolumn{1}{l|}{50} & \multicolumn{1}{l|}{} & \multicolumn{1}{l|}{50} & \multicolumn{1}{l|}{50} \\ \cline{2-3} \cline{5-6} 
\multicolumn{1}{l|}{B} & \multicolumn{1}{l|}{\textbf{100}} & \multicolumn{1}{l|}{50} & \multicolumn{1}{l|}{} & \multicolumn{1}{l|}{50} & \multicolumn{1}{l|}{\textbf{100}} \\ \cline{2-3} \cline{5-6} 
\end{tabular}
\caption{Data for the various figures \ref{mosaic-100}. We have changed one of the value to 100}
\label{table-100}
\end{center}
\end{table}

\begin{table}[h!]
\begin{center}
\begin{tabular}{l l l l l l}
                      &     C                & D                    &                       & C                    & D                     \\ \cline{2-3} \cline{5-6} 
\multicolumn{1}{l|}{A} & \multicolumn{1}{l|}{\textbf{200}} & \multicolumn{1}{l|}{50} & \multicolumn{1}{l|}{} & \multicolumn{1}{l|}{50} & \multicolumn{1}{l|}{\textbf{200}} \\ \cline{2-3} \cline{5-6} 
\multicolumn{1}{l|}{B} & \multicolumn{1}{l|}{50} & \multicolumn{1}{l|}{50} & \multicolumn{1}{l|}{} & \multicolumn{1}{l|}{50} & \multicolumn{1}{l|}{50} \\ \cline{2-3} \cline{5-6} 
                      &                       &                       &                       &                       &                       \\
                      &     C                & D                    &                       & C                    & D                     \\ 
\cline{2-3} \cline{5-6} 
\multicolumn{1}{l|}{A} & \multicolumn{1}{l|}{50} & \multicolumn{1}{l|}{50} & \multicolumn{1}{l|}{} & \multicolumn{1}{l|}{50} & \multicolumn{1}{l|}{50} \\ \cline{2-3} \cline{5-6} 
\multicolumn{1}{l|}{B} & \multicolumn{1}{l|}{\textbf{200}} & \multicolumn{1}{l|}{50} & \multicolumn{1}{l|}{} & \multicolumn{1}{l|}{50} & \multicolumn{1}{l|}{\textbf{200}} \\ \cline{2-3} \cline{5-6} 
\end{tabular}
\caption{Data for the various figures \ref{mosaic-200}. We have changed one of the value to 200}
\label{table-200}
\end{center}
\end{table}

\clearpage

\section{Conclusion}
In this paper, we have revisited the Pearson Chi-squared independence test. We have provided seven proofs for this seminal test. We also present 
an innovative and intuitive graphical representation of this test with a confidence interval. 
This enables not only interpreting visually the test but also measuring how close or far we are from accepting or rejecting the null hypothesis of non independence.
Further work could be to extend these confidence interval interpretation to contingency tables larger than two by two.

\bibliographystyle{jfe}
\bibliography{mybib2}

\begin{thebibliography}{12}
\bibitem[{Pearson(1900)}]{Pearson_1900}
Pearson, K., 1900. On the criterion that a given system of deviations from the
  probable in the case of a correlated system of variables is such that it can
  be reasonably supposed to have arisen from random sampling. Philosophical
  Magazine 302, 157--175.
\end{thebibliography}

\clearpage

\end{document}